\documentclass[11pt]{article}
\usepackage{amsfonts}
\usepackage{latexsym}
\usepackage{amssymb}
\usepackage{amscd}
\usepackage{amsmath}

\newlength{\dinwidth}
\newlength{\dinmargin}
\def\be{\begin{equation}}
\def\ee{\end{equation}}
\def\ben{\begin{displaymath}}
\def\een{\end{displaymath}}
\def\baa{\begin{eqnarray}}
\def\eaa{\end{eqnarray}}

\def\ba{\begin{array}}
\def\ea{\end{array}}

\makeatletter
\@addtoreset{equation}{section}
\makeatother

\def\phi{\varphi}

\def\a{\alpha}
\def\g{{\bf m}}

\def\b{\beta}

\def\l{\lambda}

\def\th{\vartheta}

\def\o{\omega}
\def\O{\Omega}

\def\Fcal{{\cal F}}

\def\B{{\bf \sigma}}

\def\CP1{{\mathbb C}P^1}

\def\la{\label}

\def\f{\frac}
\def\L{{\cal L}}
\def\p{\partial}

\def\log{\ln}

\def\la{\label}

\def\f{\frac}
\def\L{{\cal L}}
\def\p{\partial}
\def\res{{\rm res}}

\def\be{\begin{equation}}
\def\ee{\end{equation}}
\def\bem{\begin{equation*}}
\def\eem{\end{equation*}}
\def\ben{\begin{displaymath}}
\def\een{\end{displaymath}}
\def\baa{\begin{eqnarray}}
\def\eaa{\end{eqnarray}}

\def\ba{\begin{array}}
\def\ea{\end{array}}

\def\CP1{{\bf CP}^1}

\def\Wt{\tilde B}

\def\a{\alpha}

\def\b{\beta}

\def\e{\eta}

\def\l{\lambda}
\def\lt{\tilde \l}

\def\th{\vartheta}
\def\tht{\tilde \th}

\def\o{\omega}
\def\O{\Omega}

\def\f{\frac}
\def\la{\label}

\def\L{{\cal L}}

\def\p{\partial}

\def\Qc{{\cal Q}}

\def\res{{\rm res}}

\def\et{{\tilde{\e}}}

\def\Lt{{\tilde{\L}}}
\def\gt{{\tilde{g}}}

\newtheorem{definition}{Definition}
\newtheorem{remark}{Remark}
\newtheorem{theorem}{Theorem}
\newtheorem{corollary}{Corollary}
\newtheorem{lemma}{Lemma}
\newtheorem{proposition}{Proposition}

\begin{document}

\title{Genus one polyhedral surfaces, spaces of quadratic differentials on tori and determinants of Laplacians}
\author{ Yulia Klochko\footnote{e-mail: yulia@mathstat.concordia.ca}\,\, and Alexey
Kokotov\footnote{e-mail: alexey@mathstat.concordia.ca} }

\maketitle

\begin{center}
Department of Mathematics and Statistics, Concordia
University
7141 Sherbrooke West, Montreal H4B 1R6, Quebec,  Canada
\end{center}

\vskip0.5cm
{\bf Abstract.}
 We prove a formula for the determinant of Laplacian on an arbitrary compact polyhedral surface of genus one.
 This formula generalizes the well-known Ray-Singer result for a flat torus.
 A special case of flat conical metrics given by the modulus of a meromorphic quadratic differential
 on an elliptic surface is also considered. We study the determinant of Laplacian as a functional on the
 moduli space ${\cal}Q_1(1, \dots, 1, [-1]^L)$ of meromorphic quadratic differentials with $L$ simple poles and $L$
 simple zeros and derive formulas for variations of this functional with respect to natural
 coordinates on ${\cal}Q_1(1, \dots, 1, [-1]^L)$. We give also a new proof of Troyanov's theorem stating the existence of a conformal flat conical metric on a compact Riemann surface of arbitrary genus with a prescribed divisor of conical points.

\section{Introduction} There exist several equivalent ways to look at  compact Riemann surfaces: for instance, one can define them via  algebraic equations or make use of one of the uniformization theorems, introducing the surface as, say, the quotient of  the upper half-plane over the action of a Fuchsian group. Another possibility to get a Riemann surface comes from Riemannian geometry: a two-dimensional  Riemannian manifold carries the natural complex structure  defined via isothermal local parameters.

Another, simple and elementary, way to represent a Riemann surface is the
following:  one can consider the boundary of an arbitrary (connected but,
generally, not simply connected) polyhedron in the three dimensional
Euclidean space. This is a polyhedral surface which carries the structure of
a complex  manifold (the corresponding system of  holomorphic local
parameters is obvious for all points except the vertices; near a vertex one
should introduce the local parameter $\zeta=z^{2\pi/\alpha}$, where $\alpha$
is the sum of the angles adjacent to the vertex). In this way the Riemann
surface comes together with a conformal metric; this metric is flat and has
conical singularities at the vertices. Actually, to perform this construction
it is not necessary to start from polyhedra embedded in the three dimensional
Euclidean space, one can use instead some simplicial complex, thinking of a
polyhedral surface as glued from plane triangles.

Troyanov
(see \cite{Tro}) proved that on any compact Riemann surface there exists a flat conformal conical metric with a prescribed divisor of conical points (see the precise formulation of this theorem below). Moreover, he noticed that any compact Riemann surface with flat conformal conical metric admits a proper triangulation (i. e. each conical point is a vertex of some triangle of the triangulation). This means that the above construction is universal: any compact Riemann surface can be glued from triangles.

The goal of this paper is to study the determinant of the Laplacian (acting in the trivial line bundle over the surface) as a functional on the space of Riemann surfaces with conformal flat conical metrics (polyhedral surfaces). The similar question for {\it smooth} conformal metrics and arbitrary holomorphic bundles was very popular in the eighties and early nineties being motivated by the string theory. Among the most notable results one can mention the Ray-Singer calculation of the determinant of the Laplacian in arbitrary flat line bundle over flat tori \cite{Ray}, an explicit formula for the determinant of Laplacian in the Arakelov metric found by Dugan and Sonoda \cite{DS}, the D'Hoker-Phong formula relating the determinant of the Laplacian in the Poincar\'e metric to Selberg's zeta-function \cite{DP1}, the
Zograf-Takhtajan formula for variation of the determinant of Laplacian in the Poincar\'e metric with respect to moduli of the Riemann surface \cite{Uspehi}, Fay's formula for variation of the determinant of Laplacian under arbitrary (not necessarily conformal) variation of the metric \cite{Fay92}.

The determinants of Laplacians in flat singular metrics  are much less studied: among the very few appropriate references we mention \cite{DP}, where the determinant of the Laplacian in conical metric was defined via some special regularizations of the diverging Liouville integral and the question about the relation of such a definition with the spectrum of the Laplacian remained open, and two papers \cite{King}, \cite{Au} dealing with flat conical metrics on the Riemann sphere.

In \cite{Leipzig} the determinant of the Laplacian was studied as a functional
$${\cal H}_g(k_1, \dots, k_M)\ni (\L, \omega)\mapsto {\rm det}\,\Delta^{|\omega|^2}$$
on the space ${\cal H}_g(k_1, \dots, k_M)$ of equivalence classes of pairs $(\L, \omega)$, where
$\L$ is a compact Riemann surface of genus $g$ and $\omega$ is a holomorphic one-form (an Abelian differential) with
$M$ zeros of multiplicities $k_1, \dots, k_M$. Here ${\rm det}\,\Delta^{|\omega|^2}$ stands for the determinant of the Laplacian in the flat metric $|\omega|^2$ having conical singularities at the zeros of $\omega$.
The corresponding results for the moduli spaces $Q_g(k_1, \dots, k_M, [-1]^L)$ of quadratic  differentials with $M$ zeros of multiplicities $k_1, \dots, k_M$ and $L$ simple poles were stated in \cite{Leipzig} without proofs.
The flat conical metric $|\omega|^2$ considered in \cite{Leipzig} is very special: the divisor of the conical points of this metric is not arbitrary (it should be the canonical one, i. e. coincide with the divisor of a holomorphic one-form) and
the conical angles at the conical points are integer multiples of $2\pi$.

In the present paper we study determinants of Laplacians on arbitrary polyhedral surfaces of genus one.
Our first main result is formula (\ref{fMT}) giving an explicit expression for the determinant of the Laplacian on arbitrary polyhedral torus. Then we consider an important special
case of  flat metrics  given as the modulus of a meromorphic quadratic differential on the torus with at most simple poles.
 In this case we give  simple and straightforward proofs of the results announced in \cite{Leipzig},
 in particular, we  derive formulas of the Rauch type for variations of basic holomorphic differential and
 the period of the elliptic surface under variation of the natural holomorphic coordinates on the
 moduli space of meromorphic quadratic differentials. The second main result of the paper is
 Theorem \ref{M3} below which gives variational formulas for the determinant of the Laplacian as
 a
 functional on this moduli space.

Although in this paper we deal with elliptic surfaces only, we start it with a new proof of Troyanov's existence theorem
for flat conical metrics on Riemann surfaces of an arbitrary genus; in contrast to previously known proofs of this theorem our proof is constructive.

The first author thanks Dmitry Korotkin who taught her the subject and was giving numerous advices at different stages of  preparation of this paper, the second author thanks Max-Plank-Institute f\"ur Mathematik in den Naturwissenschaften (Leipzig) for hospitality and excellent working conditions. We both acknowledge useful discussions with Marco Bertola and Alexandre Bobenko.

\vskip0.8cm
\section{Flat conical metrics on surfaces}
\subsection{Troyanov's theorem}
Let $\sum_{k=1}^N\beta_kP_k$ be a (generalized, i. e. the coefficients $\beta_k$ are not necessary integers) divisor on a compact Riemann surface $\L$ of genus $g$. Let also $\sum_{k=1}^N\beta_k=2g-2$. Then, according to Troyanov's theorem (see \cite{Tro}), there exists a (unique up to a homothety) conformal flat metric $\g$ on $\L$
which is smooth in $\L\setminus\{P_1, \dots, P_N\}$ and has simple singularities of order $\beta_k$ at $P_k$.
The latter means that in a vicinity of $P_k$ the metric $\g$ can be represented in the form
\begin{equation}\label{met}\g=e^{u(z, \bar z)}|z|^{2\beta_k}|dz|^2,\end{equation}
where $z$ is a conformal coordinate and $u$ is a smooth real-valued function. In particular, if
$\beta_k>-1$ the point $P_k$ is conical with conical angle $2\pi (\beta_k+1)$.
Here we construct the metric $\g$ explicitly, giving an effective proof of Troyanov's theorem.

Fix a canonical basis of cycles on $\L$ (we assume that $g\geq 1$, the case $g=0$ is trivial) and let $E(P, Q)$ be the prime-form (see \cite{Fay73}).
Then for any divisor ${\cal D}=r_1Q_1+\dots r_mQ_M-s_1R_1-\dots -s_NR_N$ of degree zero on $\L$ (here
the coefficients $r_k, s_k$ are positive integers)
the meromorphic differential
$$\omega_{{\cal D}}=d_z\log\frac{\prod_{k=1}^ME^{r_k}(z,Q_k)}{\prod_{k=1}^NE^{s_k}(z, R_k)}$$
is holomorphic outside ${\cal D}$ and has the first order poles at the points of ${\cal D}$ with residues $r_k$ at $Q_k$ and $-s_k$ at $R_k$.
Since the prime-form is single-valued along the $a$-cycles, all the $a$-periods of the differential $\omega_{\cal D}$ vanish.

Let $\{v_\alpha\}_{\alpha=1}^g$ be the basis of holomorphic normalized differentials and
${\mathbb B}$ the corresponding matrix of $b$-periods.
Then all the $a$- and $b$-periods of the meromorphic differential
$$\Omega_{\cal D}=\omega_{\cal D}-2\pi i\sum_{\a, \b=1}^g ((\Im{\mathbb B})^{-1})_{\a\b}\Im\left( \int_{s_1R_1+\dots s_NR_N}^{r_1Q_1+\dots r_M Q_M}v_\beta\right)v_\alpha$$
are purely imaginary (see \cite{Fay73}, p. 4).

Obviously, the differentials $\omega_{\cal D}$ and $\Omega_{\cal D}$ have the same structure of poles: their difference is a holomorphic $1$-form.

Choose a base-point $P_0$ on $\L$ and introduce the following quantity
$${\cal F}_{\cal D}(P)=\exp\int_{P_0}^P\Omega_{\cal D}.$$
Clearly, ${\cal F}_{\cal D}$ is a meromorphic section of some {\it unitary} flat line bundle over $\L$, the divisor of this section coincides with ${\cal D}$.

Now we are ready to construct the metric $\g$. Choose any holomorphic differential $w$ on $\L$ with, say, only simple zeros $S_1, \dots, S_{2g-2}$.
Then one can set $\g=|u|^2$, where
\begin{equation}\label{reshen}u(P)=w(P)\Fcal_{(2g-2)S_0-S_1-\dots S_{2g-2}}(P)\prod_{k=1}^N\left[\Fcal_{P_k-S_0}(P)\right]^{\beta_k}\end{equation}
and $S_0$ is an arbitrary point.

Notice that in case $g=1$ the second factor in (\ref{reshen}) is absent and the remaining part is nonsingular at the point $S_0$.

\subsection{Distinguished local parameter}
In a vicinity of a conical point the flat metric (\ref{met}) takes the form
$$\g=|g(z)|^2|z|^{2\beta}|dz|^2$$
with some holomorphic function $g$ such that $g(0)\neq0$. It is easy to show (see, e. g., \cite{Tro},  Proposition 2) that there exists a holomorphic change of variable $z=z(w)$ such that in the local parameter $w$
$$\g=|w|^{2\beta}|dw|^2\,.$$
We shall call the parameter $w$ (unique up to a constant factor $c$, $|c|=1$) {\it distinguished}. In case $\beta>-1$ the existence of the distinguished parameter means that in a vicinity of conical point the surface $\L$ is isometric to the standard cone with conical angle $2\pi(\beta+1)$.

\section{Flat conical metrics on tori and determinants of Laplacians}
\subsection{Determinants of Laplacians}
From now on $\L$ is an elliptic ($g=1$) Riemann surface and it is assumed that $\L$ is the quotient of the complex plane ${\mathbb C}$ by the lattice generated by $1$ and $\sigma$, where $\Im\sigma>0$. The differential $dz$ on ${\mathbb C}$ gives rise to a holomorphic differential $v_0$ on $\L$ with periods $1$ and $\sigma$.

Let $\sum_{k=1}^N\beta_kP_k$ be a generalized divisor on $\L$ with
$\sum_{k=1}^N\beta_k=0$ and assume that $\beta_k>-1$ for all $k$. Let $\g$ be
a flat conical metric corresponding to this divisor via Troyanov's theorem.
Clearly, it has a finite area and is defined uniquely when this area is
fixed. Fixing numbers $\b_1, \dots, \b_N>-1$ such that
$\sum_{k=1}^N\beta_k=0$,  we define the space ${\cal M}(\beta_1, \dots,
\beta_N)$ as the moduli space of pairs $(\L, \g)$, where $\L$ is an elliptic
surface and $\g$ is a flat conformal metric on $\L$ having $N$ conical
singularities with conical angles $2\pi (\b_k+1)$, $k=1, \dots, N$. The space
${\cal M}(\beta_1, \dots, \beta_N)$ is a connected  orbifold of the real
dimension $2N+3$.

Let $z=x+iy$ be a conformal coordinate on $\L$ and let $\g=\rho^{-2}(z, \bar z)\widehat{dz}=\rho^{-2}dx\,dy$. Denote by $\Delta^\g$ the Friedrichs extension of the operator
$$C_0^\infty(\L\setminus\{P_1, \dots, P_N\})\ni f\mapsto 4\rho^2\partial^2_{z\bar z}f.$$

The determinant of $\Delta^\g$ for flat metrics with conical singularities was first defined in \cite{King}.  Briefly, this definition looks as follows.  Cheeger's theorem (\cite{Cheeger}) states that
the spectrum, $\{\l_k\}$, of $\Delta^\g$ is discrete (with each eigenvalue having finite multiplicity) and its counting function, $N(\l)$, obeys the standard spectral asymptotics $N(\l)=O(|\l|)$ at the infinity. Moreover, from the results of Br\"uning and Seeley \cite{BrSeel} it follows that the analytic continuation of the corresponding operator zeta-function $$\zeta_{\Delta^
\g}(s)=\sum_{\l_k\neq 0}\l_k^{-s}$$ (the latter series converges to a holomorphic function of $s$ in the half-plane $\{\Re s>1\}$) is meromorphic in the complex plane and has no pole at $s=0$.
Therefore, one can define the determinant of the operator $\Delta^\g$ via the standard Ray-Singer regularization:
$${\rm det}\Delta^\g=\exp\{-\zeta'_{\Delta^\g}(0)\}\,.$$

The main result of the present paper, stated below as Theorem 1, is an explicit formula for the function
$${\cal M}(\beta_1, \dots, \beta_N)\ni (\L, \g)\mapsto {\rm det}\Delta^\g\, .$$

Write the  normalized holomorphic differential $v_0$ on the elliptic surface $\L$ in the
the distinguished local parameter $w_k$ near the conical point $P_k$ ($k=1, \dots, N$) as
$$v_0=f_k(w_k)dw_k$$
and define
\begin{equation}{\bf f}_k:=f_k(w_k)|_{w_k=0}, \ k=1, \dots, N\,.   \end{equation}

 \begin{theorem}\label{MT}
The following formula holds true
\begin{equation}\label{fMT}
{\rm det}\Delta^\g=C|\Im \sigma|\,{\rm Area}(\L, \g)\, |\eta(\sigma)|^4\prod_{k=1}^N|{\bf f}_k|^{-\beta_k/6},
\end{equation}
where $C$ is a constant depending only on $\beta_1, \dots, \beta_N$, and $\eta$ is the Dedekind eta-function.
\end{theorem}

The proof of this theorem will be given in the next section.
\begin{remark}{\rm An analogous statement for genus $0$ polyhedral surfaces was obtained in \cite{Au}.
When the flat metric $\g$ is everywhere nonsingular formula (\ref{fMT}) reduces to the well-known Ray-Singer result \cite{Ray}.
}
\end{remark}

\subsection{Proof of Theorem \ref{MT}}
The proof uses three basic technical tools: the
Burghelea-Friedlander-Kappeler analytic surgery, the Polyakov formula and the Ray-Singer calculation of the determinant of Laplacian corresponding to smooth flat metric on the elliptic surface.

\subsubsection{Analytic surgery}
Take $\epsilon>0$ and introduce the disks $D_k(\epsilon)=\{|w_k|\leq
\epsilon\}$, centered at the conical points $P_k$, $k=1, \dots, N$. Let
$\Sigma_\epsilon=\L\setminus\cup_{k=1}^ND_k(\epsilon)$. Let also
$g_k:\overline{{\mathbb R}}_+\rightarrow {\mathbb R}$, $k=1, \dots, N$  be
smooth positive functions such that
\begin{enumerate}
\item
$\int_0^1g_k^2(r)rdr=\int_0^1r^{2\beta_k+1}dr=\frac{1}{2\beta_k+2}$,
\item
$g_k(r)=r^{\beta_k}$ for $r\geq 1$.
\end{enumerate}

Define the family of {\it smooth} conformal metrics $\g_\epsilon$ on $\L$ via
$$\g_\epsilon(z)=\begin{cases}
\epsilon^{2\beta_k}g_k^2(|w_k|/\epsilon)|dw_k|^2,\ \ \ \ \ z\in D_k(\epsilon), \ \ k=1, \dots, N\\
\g(z), \ \ \ \ \ \ \ \ z\in \Sigma_\epsilon
\end{cases}$$
The metrics $\g_\epsilon$ converge to $\g$ in $\L\setminus\{P_1, \dots, P_N\}$ as $\epsilon \to 0$  and
$${\rm Area}(\L, \g_\epsilon)={\rm Area}(\L, \g).$$

\begin{lemma}\label{SOV}
 Let $\partial_t$ be the differentiation with respect to one of the
coordinates on ${\cal M}(\beta_1, \dots, \beta_N)$
and let ${\rm
det}\Delta^{\g_\epsilon}$ be the standard $\zeta$-regularized determinant of
the Laplacian corresponding to the smooth metric $\g_\epsilon$.
Then
\begin{equation}\label{sovpad}
\partial_t\log{\rm  det}\Delta^{\g}=\partial_t\log{\rm
det}\Delta^{\g_\epsilon}.
\end{equation}
\end{lemma}
{\bf Proof.}
For simplicity suppose first that $N=1$.
Let
$(\Delta^{\g_\epsilon}|D)$ and $(\Delta^{\g_\epsilon}|\Sigma)$ be the
operators of the Dirichlet boundary problem for $\Delta^{\g_\epsilon}$ in
domains
$D:=D_1(\epsilon)$ and $\Sigma:=\Sigma_\epsilon$ respectively.
Define the Neumann jump operator (a pseudodifferential operator on
$\partial D$ of order $1$)
$R:C^{\infty}(\partial D)\to C^{\infty}(\partial D)$
by
$$R(f)=\partial_\nu (V^{-}-V^{+}),$$
where $\nu$ is the outward normal to $\partial D$, the functions
$V^{-}$ and $V^{+}$ are
the solutions of the boundary value problems
$\Delta^{\g_\epsilon} V^{-}=0$ in $D$, $V^{-}|_{\partial D}=f$ and
$\Delta^{\g_\epsilon} V^{+}=0$ in $\Sigma$, $V^{+}|_{\partial D}=f$.

In what follows it is crucial that the Neumann jump operator does not
change if we vary the metric  within the same
conformal class.
Due to Theorem $B^{*}$ from \cite{BFK}, we have
\begin{equation}\label{s1}
{\rm  det}\Delta^{\g_\epsilon}={\rm
det}(\Delta^{\g_\epsilon}|D)\,{\rm det}(\Delta^{\g_\epsilon}|\Sigma)\,{\rm
det}R\,\{{\rm Area}(\L,\g_\epsilon)\}\,\{l(\partial D)\}^{-1},
\end{equation}
where $l(\partial D)$ is the length of the contour $\partial D$ in the
metric $\g_\epsilon$
\footnote{ We have excluded the zero modes of an operator from the definition of its determinant, so we are using  the same notation ${\rm det}\,A$ for the determinants of operators $A$ with and without zero modes. In \cite{BFK} the determinant of an operator $A$ with zero modes is always equal to zero, and what we  call here ${\rm det}\, A$ in \cite{BFK} is called the modified determinant and denoted by ${\rm det}^* \,A$.  }.

Analogous statement holds if the metric defining the Laplacian has a conical
singularity inside $D$ (see \cite{Leipzig}). One has the surgery formula for
the operator $\Delta^{\g}$:
\begin{equation}\label{s2}
{\rm  det}\Delta^{\g}={\rm
det}(\Delta^{\g}|D)\,{\rm
det}\,(\Delta^{\g}|\Sigma)\,{\rm  det}R\,\{{\rm
Area}(\L, \g)\}\,\{
l(\partial D)\}^{-1}.
\end{equation}

Notice that the variations of the logarithms of the  first factors in right hand sides of
(\ref{s1}) and (\ref{s2}) vanish (these factors are independent of $t$)
whereas the variations of logarithms of all the remaining factors coincide. This
leads to (\ref{sovpad}).
To consider the general case ($N>1$) one should apply an obvious
generalization of the surgery formula for several non-overlapping discs;
similar result can be found in (\cite{Lee}, remark on page 326).
$\Box$

\subsubsection{Polyakov's formula} We state this result in the form given in (\cite{Fay92}, p. 62).
 Let $\g_0=\rho_0^{-2}(z, \bar z)\widehat{dz}$ and $\g_1=\rho_1^{-2}(z, \bar z)\widehat{dz}$ be two
{\it smooth} conformal metrics on $\L$ and let ${\rm det}\Delta^{\g_0}$ and ${\rm det}\Delta^{\g_0}$ be the determinants of the corresponding Laplacians (defined via the standard Ray-Singer regularization). Then
\begin{equation}\label{Polyakov}
\frac{{\rm det}\Delta^{\g_1}}{{\rm det}\Delta^{\g_0}}=\frac{{\rm Area}(\L, \g_1)}{{\rm Area}(\L, \g_0)}
\exp\left\{\frac{1}{3\pi}\int_\L\log\frac{\rho_1}{\rho_0}\partial^2_{z\bar z}\log(\rho_1\rho_0)\widehat{dz}\right\}\, .
\end{equation}
\subsubsection{Ray-Singer formula} Let $\Delta$ be the Laplacian on $\L$ corresponding to the flat smooth metric $|v_0|^2$, where $v_0$ is the normalized holomorphic differential. The following formula for ${\rm det}\Delta$ was proved in \cite{Ray}:
\begin{equation}\label{RS}
{\rm det}\Delta=C|\Im \sigma|^2|\eta(\sigma)|^4,
\end{equation}
where $C$ is a $\sigma$-independent constant.

\subsubsection{Proof of Theorem \ref{MT}}
By virtue of Lemma \ref{SOV} one has the relation
\begin{equation}\label{e1}
\partial_t\left\{  \log\frac{{\rm det}\Delta^\g}{{\rm Area}(\L, \g)}         -\log\frac{{\rm det}\Delta}{\Im\sigma}\right\}=
\partial_t\left\{ \log\frac{{\rm det}\Delta^{\g_\epsilon}}{{\rm Area}(\L, \g_\epsilon)}             -\log\frac{{\rm det}\Delta}{\Im\sigma}\right\}\,.
\end{equation}
Applying to the r. h. s. of (\ref{e1}) Polyakov's formula, we get
\begin{equation}
\partial_t\left\{\log\frac{{\rm det}\Delta^\g}{{\rm Area}(\L, \g)}- \log\frac{{\rm det}\Delta}{\Im\sigma}  \right\}=\sum_{k=1}^N\frac{1}{3\pi}\partial_t\int_{D_k(\epsilon)}(\log G_k)_{w_k\bar w_k}\log|f_k|\widehat{dw_k},
\end{equation}
where $G_k(w_k)=\epsilon^{-\beta_k}g_k^{-1}(|w_k|/\epsilon)$. Notice that the function $G_k$ coincides with $|w_k|^{-\beta_k}$ in a vicinity of the circle $\{|w_k|=\epsilon\}$ and the Green formula implies that
$$\int_{D_k(\epsilon)}(\log G_k)_{w_k\bar w_k}\log|f_k|\widehat{dw_k}=\frac{i}{2}\left\{
\oint_{|w_k|=\epsilon}(\log|w_k|^{-\beta_k})_{\bar w_k}\log|f_k|d\bar w_k+\right.$$
$$\left.+\oint_{|w_k|=\epsilon}\log|w_k|^{-\beta_k}(\log|f_k|)_{w_k}dw_k+\int_{D_k(\epsilon)}(\log|f_k|)_{w_k\bar w_k}\log G_kdw_k\wedge d\bar w_k
\right\}$$
and, therefore,
\begin{equation}\label{e2}\partial_t\int_{D_k(\epsilon)}(\log G_k)_{w_k\bar w_k}\log|f_k|\widehat{dw_k}=-\frac{\beta_k\pi}{2}\partial_t\log|{\bf f}_k|+o(1)\end{equation}
as $\epsilon\to 0$. Formula (\ref{fMT}) follows from (\ref{e1}), (\ref{e2}) and (\ref{RS}). $\square$

\section{Spaces of meromorphic quadratic differentials on elliptic surfaces}
Here we study reductions of formula (\ref{fMT}) to the case of flat conical
metrics $|W|$, where $W$ is a meromorphic quadratic differential on $\L$
having only simple poles.  We assume that the zeroes of $W$ are also simple,
although with a little more effort one can consider the general case of
arbitrary multiplicities. Notice that the  metric $|W|$ is flat and has
conical points with conical angles $3\pi$ at the zeroes of $W$ and $\pi$ at
the poles of $W$ and, of course, the divisor of conical points is not
arbitrary --- it should be linearly equivalent to zero (since the canonical
divisor of an elliptic surface coincides with the principle one).

Following \cite{Lanneau}, \cite{ZorMas}, introduce the space   $\Qc_1(1,...,1,[-1]^L)$ of equivalence classes of pairs
$(\L, W)$, where $\L$ is an elliptic surface and $W$ is a meromorphic quadratic  differential on $\L$ with $L$ simple zeroes and $L$ simple poles\footnote{Two pairs $(\L_1, W_1)$ and $(\L_2, W_2)$ are called equivalent if there exists a biholomorphic map $f:\L_1\to\L_2$ such that $f_*W_2=W_1$ }. The space $\Qc_1(1,...,1,[-1]^L)$ is known to be a connected complex orbifold \cite{Lanneau}. (It should be noted that the space $\Qc_1(1, -1)$ is empty.)

Notice that due to modular properties of Dedekind's eta-function the product $|\Im \sigma||\eta(\sigma)|^4$ depends only on the conformal class of the elliptic surface $\L$ (and not on the choice of the canonical basis of cycles on $\L$). So one can introduce the function
$${\cal T} : \Qc_1(1,...,1,[-1]^L)\ni (\L, W)\mapsto \frac{{\rm det}\Delta^{|W|}}{|\Im \sigma||\eta(\sigma)|^4\,{\rm Area}\, (\L, |W|)}$$
and by (\ref{fMT}) we have
$${\cal T}(\L, W)=C\,|\tau|^2,$$
with $C$ being a constant independent of $(\L, W)$ and $\tau$ given by
\begin{equation}\label{tW}\tau=\left(\frac{\prod_{k=1}^L{\bf h}_k}{\prod_{k=1}^L{\bf f}_k}\right)^{\frac{1}{24}}.\end{equation}
Here ${\bf f}_k$ (respectively ${\bf h}_k$) is the value of some chosen (say, normalized differential $v_0$)  holomorphic differential on $\L$ at the $k$-th zero
(respectively $k$-th pole) of the quadratic differential $W$ calculated in the distinguished local parameter.
Now, in contrast to Theorem 1, we split the conical points into two types  (with angle $\pi$ and with angle $3\pi$), that is why we use the new notation for the  values of $v_0$ at the conical points with angle $\pi$.

The main goal of the remaining part of this paper is to study $\tau$ as a function of moduli (the holomorphic coordinates on $\Qc_1(1,...,1,[-1]^L)$).

\subsection{Local coordinates on $\Qc_1(1,...,1,[-1]^L)$ }

 For any pair $(\L, W)$ from $\Qc_1(1,...,1,[-1]^L)$ one can
construct the so-called canonical two-fold covering
$$\pi : \Lt \to \L$$
such that $\pi * W = \o^2$, where $\o$ is a holomorphic 1-differential on
$\Lt$. This covering is ramified over the poles and zeroes of $W$.

 Let $R_1, ..., R_L$ be the zeroes of a quadratic differential
$W$ and let $S_1, ..., S_L$ be its  poles. The only zeroes of the holomorphic differential $\o$ on $\Lt$ are the double zeroes at  $R_1, ..., R_L$, therefore, one has the relation  $2 \gt - 2 = 2L$
for the genus $\gt$ of the surface $\Lt$ and $\gt = L+1$.

Denote by $*$ the holomorphic involution on $\Lt$ interchanging the sheets
of the canonical covering. The differential $\o(P)$ is anti-invariant with respect to
involution $*$: \be \o(P^*) = -\o(P). \la{omega*}\ee
Here $\o(P)$ and $\o(P^*)$ stand for values of the differential $\o$ in any local parameter lifted from the base of the canonical covering.

Due to (\cite{Fay73}, p. 85), one can choose a canonical
basis of cycles
$$\{a_\a, b_\a, a_{\a'}, b_{\a'}, a_m, b_m\},\ \ \ \ \a, \a' = 1;\ \, m = 1, ..., L-1$$
on $\Lt$ such that
\begin{itemize}
\item The pair $(\pi a_\a, \pi b_\a)$ forms a canonical basis on $\L$.
\item The following invariance properties under the involution $*$ hold: \begin{equation}\label{invari1}
a_\a^*+a_{\a'} = b_\a^* + b_{\a'} = 0 \end{equation} and \begin{equation}\label{invari2} a_m^*+a_m =
b_m^* + b_m = 0.\end{equation}
\end{itemize}
\begin{remark}{\rm
The symbols denoting the basic cycles $a_\a, b_\a, a_{\a'}, b_{\a'}$ are provided with (extrinsic) indices $\alpha, \alpha'$ in order to make our notation agree with that of \cite{Fay73}, where the base of the two-fold covering may have arbitrary genus.}
\end{remark}
\begin{remark}\label{Mum}{\rm It is convenient to keep in mind the following informal representation of the canonical covering $\Lt$:
take the standard picture of a hyperelliptic covering of the Riemann sphere branched at $2L$  points
$R_1, \dots, R_L, S_1, \dots, S_L$ with the usual canonical basis of cycles (see, e. g., \cite{Mumford}, p. 76)
$\{a_m, b_m\}$, $m=1, \dots, L-1$. Then  make two holes on two different sheets (one under another). Now the sheets are two tori and in order to get a canonical basis on the obtained two-fold covering of the torus one have to add to the cycles
$\{a_m, b_m\}$, $m=1, \dots, L-1$ two pairs of cycles $\{a_\a, b_\a\}$ and $\{a_{\a'}, b_{\a'}\}$ lying one under another on different sheets of the covering (each pair forms a canonical basis on the corresponding torical sheet).}
\end{remark}

For the corresponding basis of normalized holomorphic differentials $u_{\a},
u_{\a'}, u_m$ on $\Lt$ we have as a corollary of (\ref{invari1},
\ref{invari2}): \be u_{\a}(P^*) = -u_{\a'}(P), \qquad u_m(P^*) = -u_m(P).
\la{u ot P*} \ee

According to \cite{Lanneau},  the complex dimension of the space
$\Qc_1(1^L,[-1]^L)$ is $2L$. As it is explained in (\cite{ZorMas}, \S 4.2;
see, also, \cite{Lanneau}, \S 2 ) one can choose a system of local
coordinates on this space as follows: \be A_\a := \oint_{a_\a} \o, \qquad
B_\a := \oint_{b_\a} \o, \qquad A_m := \oint_{a_m} \o, \qquad B_m :=
\oint_{b_m} \o \la{AB coordinates} \ee for $\a = 1, m = 1, ..., L-1$. (In
\cite{Lanneau} the above coordinates are called Kontsevich's cohomological
coordinates.)

In what follows we shall refer to the cycles $\{a_m, b_m\}$ and the
coordinates $A_m, B_m$ as {\it Latin} and to the cycles  $\{a_\a, b_\a\}$ and
the coordinates $A_\a, B_\a$   as {\it Greek}.

\subsection{Projective connections and canonical meromorphic bidifferential}
Having fixed a canonical basis of cycles on a Riemann surface, one can introduce the prime-form $E(P, Q)$ and the canonical meromorphic bidifferential $B(P, Q)=d_Pd_Q\log E(P, Q)$ (see \cite{Fay73}). Recall that the canonical meromorphic bidifferential $B(P, Q)$
  is singular on the diagonal $P=Q$ and has the following local behavior as
$P \to Q$:
\be B(x(P),x(Q)) = \left( \f1{(x(P)-x(Q))^2} +\f{1}{6}S_B(x(P)) +
o(1) \right) dx(P)dx(Q) \ee Here $x(P)$ is a local parameter of a point $P
\in \L$ and the term $S_B(x(P))$ is a projective connection. This projective connection is called the {\em Bergman projective
connection}.  Recall, that a projective connection $S$ is a quantity
transforming under the coordinate change $z=z(t)$ as follows:
$$S(t) = S(z) \left( \f{dz}{dt} \right)^2 + \{z,t\},$$
where
$$\{z,t\} = \f{z'''(t)z'(t)-\f32(z''(t))^2}{(z'(t))^2}$$
is the Schwarzian derivative.

In what follows we denote by $S_B$ (respectively $\tilde{S}_B$) and $B$
(respectively $\tilde{B}$) the Bergman projective connection and the
canonical meromorphic differential on the elliptic surface $\L$ (respectively
on the canonical covering $\Lt$ of genus $\gt=L+1$). The canonical basis of
cycles on $\L$ and $\Lt$ are chosen as in the previous section.

 With $\sigma$ denoting the $b$-period of the normalized holomorphic differential $v_0$  on $\L$, introduce the function $\et$ by the
equation
$$\et(\B) = \f{d}{d\B}\log \e(\B),$$ where $\e$ is the Dedekind eta-function.
Then the canonical meromorphic bidifferential on $\L$  has the following explicit expression:
\begin{equation}\label{Bergell}B(x,y)=\Big[\wp(\int_x^y v_0) - 4\pi i \et(\B)\Big]v_0(x)v_0(y),\end{equation} where
$\wp$ is the Weierstrass $\wp$-function (see \cite{Fay73}).

\subsection{Rauch type formulas on the space $\Qc_1(1^L,[-1]^L)$}

Varying the coordinates of the pair $(\L, W)$ in the space $\Qc_1(1^L,[-1]^L)$, we change the conformal class of the elliptic surface $\L$. The following two propositions describe the behavior of the normalized holomorphic differential
$v_0$ on $\L$ under variations of the coordinates.
Let, as before, $\omega$ be the holomorphic differential on $\Lt$ such that $\omega^2=W$. Then one can introduce the following local coordinate on $\Lt$ (outside the divisor $(\omega)$):
                    $$z(P) = \int_{R_1}^P \o.$$
Below in order to simplify the notation we always make the following
agreement.

{\it Under the expression $v_0(P)$  with the argument $P$ belonging to the canonical covering
one should understand the lift $\pi_*v_0$ of the one-form $v_0$ on the base $\L$ to the canonical covering $\Lt$.
The same agreement holds for the canonical meromorphic bidifferential $B(P, Q)$ on $\L$: if $P$ (or $Q$ or both $P$ and $Q$) belongs to the canonical covering one should apply the corresponding lift. }

\begin{proposition}
If $z(P)$ is kept fixed under the differentiation then the basic differential
$v_0$ on $\L$ depends on the coordinates $A_{\a}$ and $B_{\a}$ as follows

\be \f{\p v_0(P)}{\p A_{\a}} \Big|_{z(P)} = -\f1{2\pi i} \oint _{b_{\a}}
\f{v_0(Q) B(P,Q)}{\o(Q)}, \qquad \f{\p v_0(P)}{\p B_{\a}} \Big|_{z(P)} =
\f1{2\pi i} \oint _{a_{\a}} \f{v_0(Q) B(P,Q)}{\o(Q)}\,. \la{var form ABalpha}
\ee
\la{var Th ABalpha}
\end{proposition}

{\bf Proof.}
Let us prove the first formula of (\ref{var form ABalpha}). The differential
$\f{\p v_0(P)}{\p A{\a}} \Big| _{z(P)}$ has a jump on $\Lt$ only on the cycle
$b_{\a}$ and all the $a$-periods of this differential vanish. Therefore, one can restore this
differential in terms of the canonical
meromorphic differential $\Wt (P,Q)$ on $\Lt$:
$$\f{\p v_0(P)}{\p A_{\a}} \Big|_{z(P)} = \f1{2\pi i} \oint _{b_{\a}}
\f{v_0(Q) \Wt(P,Q)}{\o(Q)}$$
(cf., \cite{Zver}).  Recall that \be b_m=-b_m^*, \quad
\o(Q^*)=-\o(Q), \quad \o_{\a}(Q^*)=-\o_{\a}(Q) \la{*change}\ee
and that the canonical meromorphic differential on $\Lt$ satisfies  \be\Wt(P^*, Q^*) = \Wt(P,Q)\la{*changeW}\ee for any $P,Q \in \L$
and  is related to the meromorphic differential $B(P,Q)$ on $\L$ as
follows:
\be B(P,Q) =\Wt(P,Q) + \Wt(P,Q^*), \quad P,Q \in \L\la{W and Wt} \ee
(see \cite{Fay73}).
Therefore,
$$\oint _{b_{\a}} \f{v_0(Q) \Wt(P,Q)}{\o(Q)} = \f12\left\{
\oint _{b_{\a}} \f{v_0(Q) \Wt(P,Q)}{\o(Q)} +\oint _{b_{\a}} \f{v_0(Q)
\Wt(P,Q^*)}{\o(Q)} \right\} = \f12 \oint _{b_{\a}} \f{v_0(Q) B(P,Q)}{\o(Q)}$$

The second formula of (\ref{var form ABalpha}) can be proved in the same way.
$\Box$

Before writing  variational formulas  with  respect to remaining Latin
coordinates we have to introduce some new notation and make an agreement
about the choice of Latin cycles.

Let us specify the form of the distinguished local parameters at
the points $S_i$ and $R_k, \quad i, k=1,...,L$ and introduce the local parameters near the same points considered as points of the canonical covering.

The distinguished local parameter (on the base $\L$) near the point $R_k$ will be denoted  by $\l_k$: one has   $$\l_k = \left( \int_{R_k}^P \o
\right)^{2/3}.$$
For a neighborhood of $R_k$ on the covering $\Lt$ we define
the local parameter $\lt_k$ to be $\lt_k=\left( \int_{R_k}^P \o
\right)^{1/3}$.

The distinguished local parameter near the $S_i$ on $\L$ will be denoted by  $\th_i$: one has $$\th_i = \left( \int_{S_i}^P \o
\right)^2.$$
For a neighborhood of $R_k$ on the covering $\Lt$ we define the local
parameter $\tht_i$ to be $\tht_i=\int_{S_i}^P \o$.

Assume for definiteness that the Latin cycles are chosen in the following
way: we split the zeros and poles $R_1, \dots, R_L, S_1, \dots, S_L$ into $L$
pairs $(R_k, S_k)$, $k=1, \dots, L$ and choose the cycle $a_k$, $k=1, \dots
L-1$ encircling the pair $(R_{k+1}, S_{k+1})$; the cycle $b_k$ intersects the
cuts $[R_1, S_1]$ and $[R_{k+1}, S_{k+1}]$ (cf. \cite{Mumford}, p. 76). Under
this assumption we have the following expressions for $z(P)$ when $P$ belongs
to the divisor $(\o)$:
\begin{equation} \begin{gathered} z(S_1)= \sum_{m=1}^{L-1} \f{A_m}{2}, \quad
z(S_2)=\f{A_1-B_1}{2}, \quad z(R_2) = -\f{B_1}{2}, \\ z(S_k)=-\f{B_{k-1}}{2}
+ \sum_{j=1}^{k-2}\f{A_j}{2}, \quad z(R_k)=-\f{B_{k-1}}{2} +
\sum_{j=1}^{k-1}\f{A_j}{2}. \la{z-coordinates}\end{gathered} \end{equation}

It will be convenient to use the following agreement: if, say, $R_k$ is the point of the divisor $(\omega)$
then $v_0(R_k)$ and $v_0'(R_k)$ are the coefficients in the expansion of $v_0$ near the point $R_k$ of the {\it canonical covering}:
$$v_0(P)=(v_0(R_k)+v_0'(R_k)\tht_k+\dots)d\tht_k.$$
Analogously, for points $P$ outside the divisor $(\omega)$:
the quantities $v(P)$ and $v'(P)$ are defined via the expansion
$$v(Q)=(v(P)+v'(P)(z(Q)-z(P))+\dots)dz(Q)$$
near the point $P$ of the canonical covering.
The expressions $\omega'(P)$,  $\omega''(P)$, $B(P, R_k)$ etc. are understood in the same way.
Now we are ready to continue the list of variational formulas.

\begin{proposition}
If $z(P)$ is kept fixed under the differentiation and  the projection of the
point $P$ on the base of canonical covering lies outside  the projection of
the contour $b_m$ on the base for the first formula and outside  the
projection of $a_m$ on the base for the second one \footnote{This refers to
the picture explained in Remark \ref{Mum}. To avoid this referring, one has
to note that the cycles $a_m$ and $-a_m^*$ (as well as $b_m$ and $-b_m^*$)
 are freely homotopic and, therefore, by virtue of Theorem 2.5 \cite{Strebel}, bound a (uniquely defined) ring domain. The point $P$ should lie outside this domain.} then the basic differential $v_0$ on $\L$ depends on the
coordinates $A_m$ and $B_m$ as follows

\be \f{\p v_0(P)}{\p A_m} \Big|_{z(P)} = -\f1{4\pi i} \oint _{b_m} \f{v_0(Q)
B(P,Q)}{\o(Q)}, \qquad \f{\p v_0(P)}{\p B_m} \Big|_{z(P)} = \f1{4\pi i} \oint
_{a_m} \f{v_0(Q) B(P,Q)}{\o(Q)}. \la{var form ABm P outside} \ee If the
projection of $P$ on the base lies inside the projection of the contour $b_m$
than the variational formula for $v_0$ with respect to $A_m$ will look as
follows: \be \f{\p v_0(P)}{\p A_m} \Big|_{z(P)} = -\f1{4\pi i} \oint _{b_m}
\f{v_0(Q)B(P,Q)}{\o(Q)} + \f12 \f{v_0'(P)\o(P) - v_0(P)\o'(P)}{\o^2(P)}.
\la{var form Am P inside} \ee Similarly, if the projection of $P$ lies inside
the  projection of the contour $a_m$ then
 \be \f{\p v_0(P)}{\p B_m}
\Big|_{z(P)} = \f1{4\pi i} \oint _{a_m} \f{v_0(Q)B(P,Q)}{\o(Q)} +\f12
\f{v_0'(P)\o(P) - v_0(P)\o'(P)}{\o^2(P)}. \la{var form Bm P inside} \ee

\la{var Th ABm}
\end{proposition}

{\bf Proof.} The proof of  formulas (\ref{var form ABm P outside}) is similar
to the proof of (\ref{var form ABalpha}) in the Proposition \ref{var Th
ABalpha}. Let us prove (\ref{var form Bm P inside}). For $P$ in a
neighborhood of the point
 $R_k$ one has the expansion
\be v_0(P) = ({\bf f}_k+{\bf f}_{k,1}\l_k(P)+\dots)d\l_k(P), \quad P\to R_k \ee
Using the
relation between the local parameters $\l$ and $\lt$ we get that $d\l_k =
2\lt_k d\lt_k$. Taking into account that
 $$d \lt_k =
\f1{3}[z(P)-z(R_k)]^{-2/3} dz = \f{dz}{3\lt_k^2},$$
we rewrite $v_0(P)$ in the following way:
$$v_0(P)=\f{2}{3} \left( \f{{\bf f}_k}{\lt_k} +{\bf f}_{k,1}\lt_k^2+ \dots \right)dz\,.$$
Differentiate this equation with respect to $B_m$ and making use of the relation
$$\f{\p\lt_k(P)}{\p B_m} = \f1{3}[z(P)-z(R_k)]^{-2/3}\f{\p z(R_k)}{\p B_m}$$
and formulas (\ref{z-coordinates}), we see that the differential $\frac{\partial v_0(P)}{\partial B_m}$ has the pole
of the second order at $R_{m+1}$,
$$\f{\p v_0}{\p B_m} \Big|_{z(P)} = -\f1{3} v_0(R_{m+1})\f{d\lt_{m+1}}{\lt_{m=1}^2}+\dots,$$
and the only other singularity of $\frac{\partial v_0(P)}{\partial B_m}$ on $\Lt$ is the jump on the cycle $a_m$.
Thus,
\be\f{\p v_0}{\p B_m}\Big|_{z(P)} = \f1{2\pi i}\oint_{a_m}\f{v_0(Q)\Wt (P,Q)}
{\o(Q)} -\f1{3} {\bf f}_{m+1}\Wt(P,R_{m+1}) \la{v0 wrt B_m in terms of Wt}.\ee
Then from (\ref{*change}), (\ref{*changeW})
and (\ref{W and Wt}) it follows that:
$$\oint_{a_m} \f{v_0(Q)\Wt (P,Q^*)}{\o(Q)} = \oint_{a_m^*} \f{v_0(Q^*)\Wt (P,Q)}{\o(Q^*)}
= -\oint_{a_m^*} \f{v_0(Q)\Wt (P,Q)}{\o(Q)}.$$ Therefore,
\bem \begin{aligned} &\oint_{a_m} \f{v_0(Q)\Wt (P,Q)}{\o(Q)} - \oint_{a_m}
\f{v_0(Q)\Wt (P,Q^*)}{\o(Q)} = \oint_{a_m} \f{v_0(Q)\Wt (P,Q)}{\o(Q)} +
\oint_{a_m^*}
\f{v_0(Q)\Wt (P,Q)}{\o(Q)}  \\
& =2\pi i \left[\res \Big|_{Q=P} \f{v_0(Q)\Wt
(P,Q)}{\o(Q)} + \res \Big|_{Q=R_{m+1}} \f{v_0(Q)\Wt (P,Q)}{\o(Q)}\right]  \\
&= 2\pi i \left[\f{v_0'(P)\o(P) -v_0(P)\o'(P)}{\o^2(P)} + \f{1}{3}
v_0'(R_{m+1})\Wt(P,R_{m+1}) \right].
\end{aligned} \eem
Hence,
\be\oint_{a_m} \f{v_0(Q)\Wt (P,Q)}{\o(Q)} = \f12 \oint_{a_m} \f{v_0(Q) (\Wt
(P,Q^*)+\Wt(P,Q))}{\o(Q)}  $$ $$+ \pi i
\left[\f{v_0'(P)\o(P)-v_0(P)\o'(P)}{\o^2(P)}+
\f{1}{3}v_0'(R_{m+1})\Wt(P,R_{m+1})\right] \la{int}\,. \ee
 Finally,
substituting (\ref{int}) into (\ref{v0 wrt B_m in terms of Wt}) we arrive at
(\ref{var form Bm P inside}). Similarly, one can prove  formula (\ref{var
form Am P inside}). $\Box$

Integrating  formulas (\ref{var form ABalpha}) and (\ref{var form ABm P
outside}--\ref{var form Bm P inside}) over the $b$-cycles of $\L$, we get the
following result which presents an analog of the well-known Rauch formulas.
\begin{corollary}\label{Rauch}
The  $b$-periods $\sigma$ of the Riemann surface $\L$ depend on the
coordinates $A_\a$, $B_\a$, $A_m$, $B_m$ as follows:

$$\f{\p\sigma}{\p A_\a} = -\oint _{b_\a} \f{v_0^2}{\o}, \quad \quad
\f{\p\sigma}{\p B_\a} = \oint _{a_\a} \f{v_0^2}{\o},$$
$$\f{\p\sigma}{\p A_m} = -\f{1}{2} \oint _{b_m} \f{v_0^2}{\o}, \quad \quad
\f{\p\sigma}{\p B_m} = \f{1}{2} \oint _{a_m} \f{v_0^2}{\o}.$$
\end{corollary}

Our last technical result is the list of variational formulas for quantities ${\bf f}_k$ and ${\bf h}_k$.

\begin{lemma}\label{pomog}
The following variational formulas hold:
\be \f{\p {\bf f}_k}{\p A_{\a}} = -\f1{2 \pi i} \oint_{b_{\a}} \f{v_0(Q)
B(R_k,Q)}{\o(Q)}, \quad \quad \f{\p {\bf h}_i}{\p A_{\a}} = -\f1{2 \pi i} \oint_{b_{\a}}
\f{v_0(Q) B(S_i,Q)}{\o(Q)} \la{fk i hi po Aalpha} \ee

\be \f{\p {\bf f}_k}{\p B_{\a}} = \f1{2 \pi i} \oint_{a_{\a}} \f{v_0(Q)
B(R_k,Q)}{\o(Q)}, \quad  \quad \f{\p {\bf h}_i}{\p B_{\a}} = \f1{2 \pi i} \oint_{a_{\a}}
\f{v_0(Q) B(S_i,Q)}{\o(Q)} \la{fk i hi po Balpha} \ee

\be \f{\p {\bf f}_k}{\p A_m} = -\f1{4 \pi i} \oint_{b_m} \f{v_0(Q)
B(R_k,Q)}{\o(Q)}, \quad \quad \f{\p {\bf h}_i}{\p A_m} = -\f1{4 \pi i} \oint_{b_m}
\f{v_0(Q) B(S_i,Q)}{\o(Q)} \la{fk i hi po Am} \ee

\be \f{\p {\bf f}_k}{\p B_m} = \f1{4 \pi i} \oint_{a_m} \f{v_0(Q)
B(R_k,Q)}{\o(Q)}, \quad \quad \f{\p {\bf h}_i}{\p B_m} = \f1{4 \pi i} \oint_{a_m}
\f{v_0(Q) B(S_i,Q)}{\o(Q)} \la{fk i hi po Bm} \ee
\end{lemma}

{\bf Proof.} The proofs of these formulas are similar, let us prove, say, the
second formula of (\ref{fk i hi po Bm}). The proof splits into two cases
depending whether the projection of the point $P$ on the base of the
canonical covering lies inside or outside of the projection of the basic
cycle $a_m$. For brevity consider only the case when the projection of $P$
lies inside the projection of $a_m$.  In the neighborhood of $S_{m+1}$ one
has the expansion
$$v_0(P)=2\left[{\bf h}_{m+1}\left(z(P)+\f{B_{m}}{2}-\sum_{j=1}^{m-1}\f{A_j}{2}\right)
+\dots \right]dz.$$ Differentiating this equality with respect to $B_m$ and
using the variational formula (\ref{var form Bm P inside}) for $v_0$ we get

$$\f1{4 \pi i} \oint_{a_m} \f{v_0(Q)B(S_{m+1},Q)}{\o(Q)}d\th_{m+1}+\f12\f{v_0'(S_{m+1})\o(S_{m+1})-
v_0(S_{m+1})\o'(S_{m+1})}{\o^2(S_{m+1})}  $$ $$=2\left[
\f12{\bf h}_{m+1}+{\bf h}_{m+1\, B_m}'\left(z(P)+\f{B_m}{2}-\sum_{j=1}^{m-1}\f{A_j}{2}\right)+\dots
\right]dz, \quad P\to S_{m+1}$$

Notice that $v_0(S_{m+1})=0$ (recall that this is true on the canonical covering and not on the base, where the differential $v_0$ has neither zero nor poles) and $d\th_{m+1}$ can be rewritten in terms of
$z$-coordinate as  $d\th_{m+1}=2(z(P)-z(S_{m+1}))dz$. Hence,
$$\f1{4 \pi i} \oint_{a_m}\f{v_o(Q) B(S_{m+1},Q)}{\o(Q)}\cdot2(z(P)-z(S_{m+1}))+{\bf h}_{m+1}=
{\bf h}_{m+1}+2{\bf h}_{m+1\,B_m}'+\dots, \quad P\to S_{m+1}.$$
 Taking the limit $P \to S_{m+1}$, we obtain formula (\ref{fk i hi po
Bm}). $\Box$

\subsection{Wirtinger tau-function on $\Qc_1(1,...,1,[-1]^L)$}
Let $\xi:{\mathbb C}\longrightarrow {\mathbb C}/\{1, \sigma\}=\L$ be the
natural projection and let $x$ be some local parameter on $\L$. Then the
Schwarzian derivative $\{\xi^{-1}(x), x\}$, being independent of the choice
of the branch of the multivalued map $\xi^{-1}$, defines a projective
connection on $\L$. This projective connection is called (see, e. g.,
\cite{Tyurin}) {\it the invariant Wirtinger projective connection}: in
contrast to the Bergman projective connection it does not depend on the
choice of canonical basis of cycles on $\L$. In what follows we denote this
projective connection by $S_{{\rm Wirt}}$. One can also put into
correspondence to a quadratic differential $W$ on $\L$ a projective
connection $S_{\omega}$ on $\L$ via the equation
\begin{equation}\label{prcon1}S_{\omega}(x(P))=\left\{\int^P\omega, x(P)\right\}.\end{equation}
(The Schwarzian derivative at the r. h. s. is independent of the choice of the branch of $\omega=\sqrt{W}$.)

Notice that the difference between two projective connections
$S_{{\rm Wirt}}$ and $S_{\o}$ is a meromorphic quadratic differential on $L$ with poles
at the zeroes of $W$. This quadratic differential can be lifted to $\Lt$, so we may define the
the following quantities:
$$ H_{A_\a} = \f1{12 \pi i}\oint_{b_{\a}} \f{S_{{\rm Wirt}} - S_{\o}}{\o}, \quad
H_{B_\a} = -\f1{12 \pi i} \oint_{a_{\a}} \f{S_{{\rm Wirt}} - S_{\o}}{\o},$$
$$ H_{A_m} = \f1{24 \pi i} \oint_{b_m} \f{S_{{\rm Wirt}} - S_{\o}}{\o}, \quad
H_{B_m} = -\f1{24 \pi i} \oint_{a_m} \f{S_{{\rm Wirt}} - S_{\o}}{\o}$$

\begin{lemma}\label{forma} Introduce the 1-form  by
$$\O = H_{A_\a}dA_\a + H_{B_\a}dB_\a + \sum_{m=1}^{L-1}(H_{A_m}dA_m + H_{B_m}dB_m).$$
Then
\begin{itemize}
\item
the 1-form $\O$ is independent of the choice of the canonical basis with properties (\ref{invari1}, \ref{invari2}) and therefore is defined on the space $\Qc_1(1,...,1,[-1]^L)$.
\item $d\O=0$.
\end{itemize}
\end{lemma}
In the next section we shall prove that \begin{equation}\label{formazam}\O=d\log \tau,\end{equation} where $\tau$ is given by (\ref{tW}).
Since $\tau$ is a (multivalued) function on $\Qc_1(1,...,1,[-1]^L)$ having at most constant multiplicative twists along nontrivial loops in  $\Qc_1(1,...,1,[-1]^L)$ (actually its $24$-th power is single-valued on $\Qc_1(1,...,1,[-1]^L)$), equation (\ref{formazam}) implies the Lemma.

However, we notice that the direct proof of the Lemma is also possible: the first statement follows from a somewhat cumbersome calculation which uses nothing but linear algebra, whereas the second one can be proved via Rauch type formulas and manipulations with singular double integrals -- the proof of a similar statement can be found in \cite{Leipzig}.

From Lemma \ref{forma} it follows that the connection
$$d_{{\rm Wirt}}=d+\Omega$$
in the trivial line bundle over $\Qc_1(1,...,1,[-1]^L)$ is flat. This flat connection defines a character of the fundamental group of $\Qc_1(1,...,1,[-1]^L)$ which in its turn defines a flat line bundle $\Xi$ over $\Qc_1(1,...,1,[-1]^L)$.
\begin{definition}
A horizontal holomorphic section of the bundle $\Xi$ is called Wirtinger tau-function \footnote{It should be noted that its direct analog in case when the space of quadratic differentials on tori is replaced by the moduli space of meromorphic functions on tori  has the meaning of the isomonodromic tau-function of Jimbo-Miwa \cite{IMRN2}.} on the space $\Qc_1(1,...,1,[-1]^L)$.
\end{definition}

In the next section the Wirtinger tau-function will be identified with the (multivalued) function $\tau$ from (\ref{tW}).

\subsection{Calculation of Wirtinger tau-function.}
The following proposition gives an explicit expression for the Wirtinger tau-function on $\Qc_1(1,...,1,[-1]^L)$.
\begin{proposition}\label{www} Let a pair $(\L, W)$ belong to the space $\Qc_1(1,...,1,[-1]^L)$.
The Wirtinger tau-function on the
stratum $\Qc_1(1,...,1,[-1]^L)$ of the space of quadratic differentials over
the Riemann surface $\L$ is given by the expression \begin{equation}\label{Me}\tau(\L, W) =\left[
\f{\prod_{k=1}^L {\bf h}_k}{\prod_{i=1}^L {\bf f}_i} \right]^{1/24}.\end{equation}
In particular, the $24$-th power of $\tau$ is a single-valued holomorphic function on $\Qc_1(1,...,1,[-1]^L)$.
\end{proposition}
{\bf Proof.} Let
$$\bold T(A_\a, \{A_m\}):= \log \left\{ \f {\prod_{k=1}^L {\bf h}_k}{\prod_{i=1}^L\
{\bf f}_i} \right\} \quad(=24\log\tau).$$ Define the (multivalued) map $R:
t\mapsto z$ by $z=\int^P\omega$ and $t=\int^Pv_0$. Clearly, the derivative
$R'(t)$ is a single-valued function. Then the one-form $(S_{{\rm
Wirt}}-S_\omega)/\omega$ can be rewritten as $$-\f{\{R,t\}}{R'}dt,$$ where
$\{R, t\}$ is the Schwarzian derivative, and, therefore, the statement of the
proposition is equivalent to the following equalities:
$$\f{\p T}{\p A_{\a}} = -\f2{\pi
i}\oint_{b_{\a}}\f{\{R,t\}}{R'}dt, \quad \f{\p T}{\p B_{\a}} =\f2{\pi
i}\oint_{a_{\a}}\f{\{R,t\}}{R'}dt,$$
$$\f{\p T}{\p A_m} = -\f1{\pi i}\oint_{b_m}\f{\{R,t\}}{R'}dt, \quad
\f{\p T}{\p B_m} = \f1{\pi i}\oint_{a_m}\f{\{R,t\}}{R'}dt.$$
The proof of these four formulas coincide verbatim. For example,
let us prove the first one.

Using Lemma \ref{pomog} and the representation (\ref{Bergell}) of the canonical meromorphic bidifferential on an elliptic surface, we get
$$\f{\p T}{\p A_{\a}} = \sum_{k=1}^L \f{{\bf h}_k'}{{\bf h}_k} - \sum_{i=1}^L
\f{{\bf f}_i'}{{\bf f}_i} =\f1{2\pi i} \oint_{b_\a} \f{v_0(Q)}{\o(Q)} \Big\{
\sum_{k=1}^L \f{B(R_k,Q)}{{\bf f}_k} - \sum_{i=1}^L\f{B(S_i,Q)}{{\bf h}_i}
\Big\} = $$
$$\f1{2\pi i}\oint_{b_\a} \Big\{\sum_{k=1}^L \f{v_0(Q)}{\o(Q) {\bf f}_k d\l_k(P)}
\Big[\wp(\int_P^Q v_0) - 4\pi i \et(\B)\Big]v_0(P)v_0(Q)\Big\}\Big|_{P=R_k}
$$
$$-\f1{2\pi i}\oint_{b_\a} \Big\{\sum_{i=1}^L \f{v_0(Q)}{\o(Q) {\bf h}_i d\th_i(P)}
\Big\{\wp(\int_P^Q v_0) - 4\pi i\et(\B)\Big\}v_0(P)v_0(Q)\Big\}\Big|_{P=S_i}
=$$
$$-\f1{2\pi i}\oint_{b_\a}\f{v_0^2(Q)}{\o(Q)} \sum_{k=1}^L \Big[\wp(\int_{S_k}^Q v_0)
- \wp(\int_{R_k}^Q v_0) \Big].$$
Observe that the sum under the last integral coincides with $$\f{d}{dt} \Big( \f{{\cal R}''(t)}{{\cal R}'(t)} \Big),$$ where  ${\cal R}'$
is defined by the relation $W={\cal R}'(t)(dt)^2$.

 Since ${\cal R}'(t)=[R'(t)]^2$, we get
\begin{equation}\label{1}\f{\p T}{\p A_{\a}}= -\f1{\pi i}\oint_{b_\a} \Big(\f{R'''}{(R')^2} - \f{(R'')^2}{(R')^3}\Big)\,dt\,.\end{equation}
It remains to notice that
\be\oint_{b_\a}\f{R'''}{(R')^2}dt = -\oint_{b_\a}R''d\left(\f{1}{(R')^2}\right)=
2\oint_{b\a}\f{(R'')^2}{(R')^3}dt\, ,\la{2}\ee
\be\oint_{b_\a}\f{\{R,t\}}{R'}=\oint_{b_\a}\f{R'''}{(R')^2}dt -
\f{3}{2}\oint_{b_{\a}} \f{(R'')^2}{(R')^3}dt =
2\oint_{b_{\a}} \f{(R'')^2}{(R')^3}dt- \f{3}{2}\oint_{b_{\a}}
\f{(R'')^2}{(R')^3}dt = \f{1}{2}\oint_{b_{\a}} \f{(R'')^2}{(R')^3}dt\la{3}\ee
and  the desired statement follows. $\Box$
\subsection{Variational formulas for the determinant of the Laplacian}
Let a pair $(\L, W)$ belong to $\Qc_1(1,...,1,[-1]^L)$. Introduce the quantity
$$Q(\L, W)=\frac{{\rm det}\Delta^{|W|}}{\{\Im \sigma\} {\rm Area}(\L, |W|)}$$
(this is the inverse to the Quillen norm on the determinant line).
The following Theorem describes variations of $Q(\L, W)$ with respect to coordinates on the space $\Qc_1(1,...,1,[-1]^L)$.
\begin{theorem}\label{M3} The variational formulas hold:
$$\f{\partial \log Q} {\p A_\a} = \f1{12 \pi i}\oint_{b_{\a}} \f{S_B - S_{\o}}{\o}, \quad \quad
\f{\partial \log Q}{\p B_\a} = -\f1{12 \pi i} \oint_{a_{\a}} \f{S_B - S_{\o}}{\o},$$
$$ \f{\partial \log Q}{\p A_m} = \f1{24 \pi i} \oint_{b_m} \f{S_B - S_{\o}}{\o}, \quad \quad
\f{\partial \log Q}{\p B_m} = -\f1{24 \pi i} \oint_{a_m} \f{S_B - S_{\o}}{\o},$$
where $m=1, \dots, L-1$ and $S_B$ is the Bergman projective connection.
\end{theorem}
{\bf Proof.} Recall that there is the following relation between the invariant Wirtinger and the Bergman projective connections on the elliptic surface $\L$:
\begin{equation}\label{WB}
S_{{\rm Wirt}}(x)=S_B(x)+24\pi i \et(\sigma) v_0^2(x)
\end{equation}
(see, e. g., \cite{Fay73} p. 35;  since Fay uses another normalization of the basic differential, the coefficient near $\et v_0^2$ in (\ref{WB}) differs from that in \cite{Fay73}).
By virtue of Proposition \ref{www}, relation (\ref{WB}) and the Rauch formula from
Corollary \ref{Rauch}, we have
$$\f{\partial \log Q} {\p A_\a}=\f{\partial \log(|\eta(\sigma)|^4|\tau(\L, W)|^2)}{\p A_\a}=
\f{\partial \log (\eta^2(\sigma)\tau(\L, W))}{\p A_\a}=
 \f1{12 \pi i}\oint_{b_{\a}} \f{S_{{\rm Wirt}} - S_{\o}}{\o}+2\et(\sigma)\f{\p \sigma}{\p A_\a}=$$
$$ \f{1}{12 \pi i}\oint_{b_{\a}} \f{S_{{\rm Wirt}}-S_{\o}}{\o}-2\et(\sigma)\oint_{b_{\a}}\frac{v_0^2}{\o}= \f{1}{12 \pi i}\oint_{b_{\a}}\f{S_B - S_{\o}}{\o},$$
which gives the first variational formula. The remaining variational formulas can be proved in the same way. $\Box$

\section{Summary and outlook}
In this paper we study the determinant of the Laplacian on a polyhedral
surface of genus one. The method we use here (see the proof of Theorem 1) can
be considered as a generalization of the Polyakov formula, relating the
determinants of Laplacians in two smooth conformal metrics, to the case when
one of the metrics is flat conical and another is flat and everywhere
nonsingular.

Using a further generalization of the Polyakov formula to the case of two
flat conical metrics and the results of \cite{Leipzig}, it is possible to
write a closed expression for the determinant of Laplacian on a polyhedral
surface of an arbitrary genus. We hope to address this question in the near
future.

It is also interesting to look at extremal properties of the determinants of
Laplacians in conical metrics; the only known result in this direction is
contained in \cite{KokKor}, where it was solved the problem of the
maximization of the determinant of the Laplacian on the Riemann sphere over
the set of  flat metrics of area $1$ with four conical points of conical
angle $\pi$.

\end{document}